\documentclass[reqno]{amsart}
\usepackage{amssymb,latexsym}
\usepackage{amsmath}
\usepackage{amscd}
\usepackage{graphicx}
\usepackage{gensymb}
\usepackage{enumerate}
\usepackage{newlattice}
\usepackage{tikz-cd}
\usepackage{mathrsfs}
\usepackage{mathbbol}
\numberwithin{equation}{section}

\theoremstyle{plain}

\newtheorem{theorem}{Theorem}
\newtheorem{lemma}[theorem]{Lemma}

\theoremstyle{definition}

\theoremstyle{remark}

\DeclareMathOperator{\lc}{lc}
\DeclareMathOperator{\rc}{rc}

\newcommand{\three}{\SM 3}
\newcommand{\Csquare}{\SC 2^2}

\newcommand{\seven}{\SfS 7}
\newcommand{\re}{\tup{re}}

\newcommand{\Cll}[1]{\tup{C}_{\tup{ll}}(#1)}

\newcommand{\Cul}[1]{\tup{C}_{\tup{ul}}(#1)}
\newcommand{\Cur}[1]{\tup{C}_{\tup{ur}}(#1)}
\newcommand{\spec}[1]{{\tup{spec}(#1)}}
\newcommand{\Ug}{\uparrow\!}

\begin{document}

\title[Rectangular lattices and two convex sublattices]
{Homomorphisms of distributive lattices \\as restrictions of 
congruences. \\III. Rectangular lattices and\\ two convex sublattices}

\author[G.\ Gr\"atzer]{George Gr\"atzer}
\email{gratzer@me.com}
\urladdr{http://server.maths.umanitoba.ca/homepages/gratzer/}

\author[H. Lakser]{Harry Lakser}
\email[H. Lakser]{hlakser@gmail.com} 
\address{Department of Mathematics\\University of Manitoba\\Winnipeg, MB R3T 2N2\\Canada}


\dedicatory{To the memory of E.\,T. \lp Tomi\rp Schmidt,\\
whose ideas still inspire us}

\begin{abstract} 
Let $L$ be a finite  lattice and let $I$ be an ideal of $L$.
Then the restriction map is a bounded lattice homomorphism 
of the congruence lattice of~$L$ into the congruence lattice of $I$. 
In a 2009 paper, the authors proved the converse.

In a 2012 paper, G. Cz\'edli proved an analogous result for rectangular lattices.
In this paper, we prove a stronger form of Cz\'edli's result 
and provide a short, elementary, and direct proof.
\end{abstract}

\maketitle   

\section{Introduction}\label{S:Introduction}

Let $L$ be a finite lattice and let $K$ be a convex sublattice of $L$.
Then the restriction map~$\re$, defined as
$\re \colon \bga \to \bga \restr K$ 
(the congruence $\bga$ of $L$ is mapped to its restriction to $K$),
is a bounded (that is, $\set{0,1}$-) lattice homomorphism 
of $\Con{L}$ into $\Con{K}$. 

G. Gr\"atzer and H.~Lakser \cite{GL86} proved the converse.

\begin{theorem}\label{T:old1}
Let $D$ and $E$ be finite distributive lattices and let $\gf \colon
D \to E$ be a~bounded lattice homomorphism.  
Then there exist a finite lattice $L$, 
a convex sublattice $G$ of $L$ that can be chosen to be either
an ideal or a filter of $L$, and isomorphisms 
$\ga \colon D \iso \Con L$ and $\gb \colon E \iso \Con G$,
making the following diagram commutative\tup{:}
\[ 
   \begin{CD}
      D                 @>\ga>\iso>            \Con L\\
      @V\gf VV                              @V \re VV\\
      E                 @>\gb>\iso>            \Con G
      \end{CD}
\]
where $\re$ is the restriction map: for a congruence $\bga$ of $L$,
$\re(\bga)$ is $\bga$ restricted to~$G$.
\end{theorem}

See E.\,T. Schmidt~\cite{eS14} for an alternative proof.

Theorem~\ref{T:old1} is an \emph{abstract/abstract} result. 
The congruence lattices are given as abstract finite distributive lattices, 
the finite lattices $L$ and $G$ are constructed.

We can improve on Theorem~\ref{T:old1} in two ways.

\emph{First}, we can construct the finite lattices $L$ and $G$ in smaller classes of lattices.
G. Gr\"atzer and H. Lakser \cite{GL94}  constructed them as \emph{planar lattices},
while G. Gr\"atzer and H.~Lakser \cite{GL09} represented them as isoform lattices. 
16 years later, G. Gr\"atzer and E. Knapp \cite{GKn09a} 
represented finite distributive lattices 
as congruence lattices of rectangular lattices.
Theorem~\ref{T:old1} for \emph{rectangular lattices} 
was done in  G.~Cz\'edli \cite{gC12}.
(Using filters instead of ideals, see Section~\ref{S:embeddabilty}.)

\emph{Second}, we can try to loosen the abstract/abstract representation. 
This was first done in  G. Cz\'edli \cite{gC12}.

\begin{theorem}\label{T:old3}
Let $G$ be a rectangular lattice, 
and let $\gf$ be a $\set{0, 1}$-lattice homomorphism 
of a finite distributive lattice $D$ to $\Con G$. 
Then there is a rectangular lattice~$L$ containing $G$
as a convex sublattice, which can be chosen as a filter of~$L$,
and a lattice isomorphism $\ga \colon D \to \Con L$ 
such that $ \gf = \re \circ \ga$,
where $\re$ is the restriction map of congruences from $L$ to $G$. That is, the
diagram
\[ 
   \begin{CD}
      D                 @>\ga>\iso>            \Con L\\
      @V=VV                              @V \re VV\\
      D                 @>\gf>>            \Con G
      \end{CD}
\]
is commutative.
\end{theorem}
Thus, $G$ is made into a filter and $\gf$ is represented as restriction of congruences to $G$.
 This is a concrete/abstract result: 
 instead of the finite distributive latttice $E$, 
 we are given the rectangular lattice $G$ and $\Con G$ plays the role of $E$.

E.\,T. Schmidt~\cite{eS11} proved this result for the special case when $\gf$
is injective. 
 
 In this paper, we present the following result.
 
 \begin{theorem}\label{T:Main}
Let $F$ and $G$ be rectangular lattices 
and let $\gf \colon \Con{F} \to \Con{G}$
be a~bounded lattice homomorphism. 
Then there is a rectangular lattice~$L$ containing~$F$ and $G$
as convex sublattices, where $G$ can be chosen as a filter of $L$,
satisfying the following two conditions\co
\begin{enumeratei}
\item the lattice $L$ is a congruence-preserving extension of $F$\tup{;}
\item for $\bga$ of $\Con L$,  
the map $\gf$ is $\bga\restr F \to \bga\restr G$.
\end{enumeratei}
\end{theorem}

This is a concrete/concrete result: 
 instead of the finite distributive latttice $D$, 
 we are given the rectangular lattice $F$ and $\Con F$ plays the role of $D$;
 instead of the finite distributive latttice $E$, 
 we are given the rectangular lattice $G$ and $\Con G$ plays the role of $E$.
 The lattices constructed are rectangular and $G$ is constructed as a filter,
 so Theorem~\ref{T:Main} is stronger form of Cz\'edli's result.
 Our approach also provides a short, elementary, and direct proof
 based on the construction in G. Gr\"atzer and E.\,T. Schmidt~\cite{GS14}.

\section{Preliminaries}\label{S:Preliminaries}

\subsection{Notation}\label{S:Notation}

We use the notation as in \cite{CFL2}.
You can find the complete

\emph{Part I. A Brief Introduction to Lattices} and  
\emph{Glossary of Notation}

\noindent of \cite{CFL2} at 

\verb+tinyurl.com/lattices101+

\subsection{BRT}\label{S:BRT}

The Birkhoff Representation Theorem (BRT, for short)
relates bounded homomorphisms of finite distributive lattices 
and isotone maps of finite ordered sets.

For a finite distributive lattice $D$, 
let $\Ji D$ denote the ordered sets of join-irreducible elements of $D$;
for $a \in D$, let $\spec a$ denote the set of join-irreducible elements $\leq a$.

\begin{theorem}[BRT]\label{T:BRT}
Let $D$ and $E$ be finite distributive lattices and set $P = \Ji D$  and $Q = \Ji E$.
Then the following five statements hold.
\begin{enumeratei}
\item  With every bounded homomorphism $\gf  \colon D \to E$, we can associate an isotone map $\Ji \gf  \colon Q \to P$ defined by
\[
   \Ji \gf  x  =  \MMm{e \in D}{x \leq \gf e}
\]
for $x \in Q$.
\item  With every isotone map $\gy \colon Q \to P$, we can associate a
bounded homo\-mor\-phism $\Down(\gy) \colon D \to E$ defined by
\[
  \Down(\gy)(e) =  \JJ \gy^{-1}(\spec e)
\]
for $e \in D$.
\item The constructions of \tup{(i)} and \tup{(ii)} are
inverses to one another, and so yield together a bijection between
bounded homomorphisms $\gf  \colon D \to E$
and isotone maps  $\Ji \gf  \colon Q \to P$.
\item $\gf$ is one-to-one iff $\Ji \gf$ is onto.
\item $\gf$ is onto if{}f $\Ji \gf$ is an order-embedding.
\end{enumeratei}
\end{theorem}

The two main formulas in this result are easy to visualize. 
For the formula in (i),  
for $x \in Q$, take the set $X$ of all elements of $D$ mapped by~$\gf$ onto~$x$ or above.
Since $\gf$ is a homomorphism, 
it follows that the meet of all these elements, let us denote it by $x^\dag$, 
is still mapped by $\gf$ onto $x$ or above.
The map $ \Ji \gf$ maps $x$ into $x^\dag$.

With every isotone map $\gy \colon Q \to P$, we can associate a
bounded homo\-mor\-phism $\Down(\gy) \colon D \to E$ defined by
\[
  \Down(\gy)(e) =  \JJ \gy^{-1}(\spec e)
\]
for $e \in D$.

Now for the formula in (ii), take an element $e \in D$.
We form the set $Y$ of join-irreducible elements of $E$
mapped by $\gy$ to an element $\leq e$. 
The join of $Y$ denoted by $e^ \ddag$, is an element of $E$.
The map  $\Down(\gy)$ maps $e$ into $e^ \ddag$.

The constructions of \tup{(i)} and \tup{(ii)} are
inverses to one another, and so yield together a~bijection between
bounded homomorphisms $\gf  \colon D \to E$ and isotone maps  $\gy  \colon Q \to P$.

\subsection{Triple gluing}\label{S:Triple}

For rectangular lattices, we need the following variant of gluing.

\begin{figure}[bth]
\centerline{\includegraphics[scale=0.9]{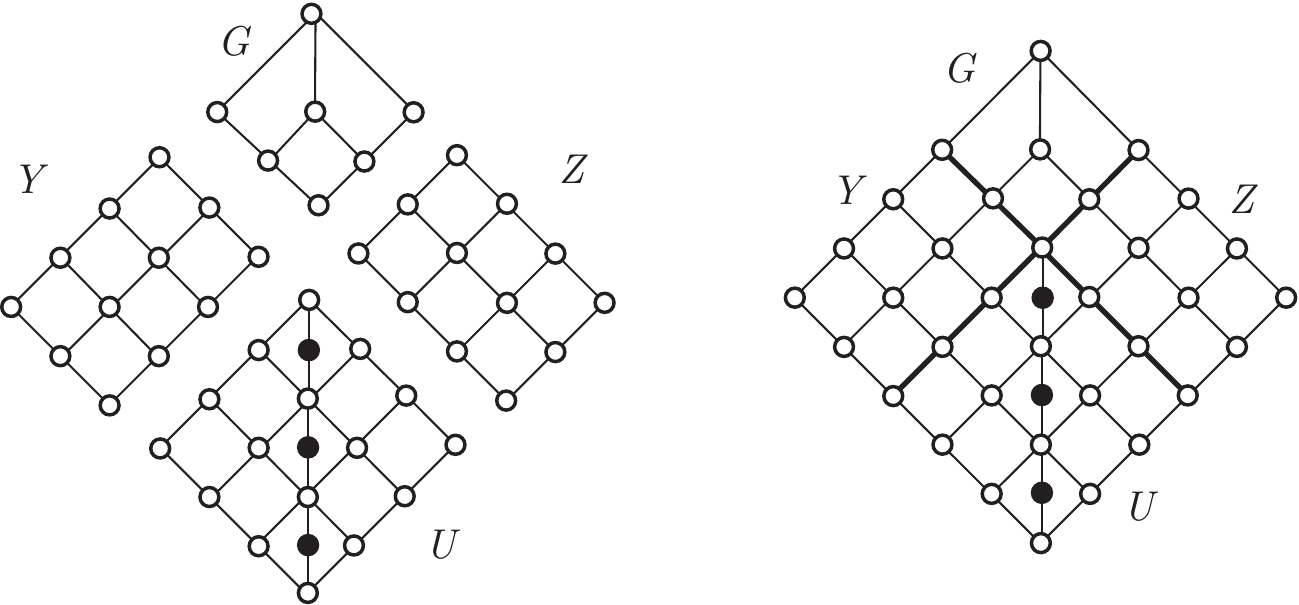}}
\caption{Triple gluing}
\label{F:first1}
\end{figure}

Let $G, Y,  Z, U$ be rectangular lattices, 
arranged in the plane as in Figure~\ref{F:first1}, 
such that the facing boundary chains have the same number of elements. 
First, we glue $U$ and $Y$ together over the boundaries, 
to~obtain the rectangular lattice $X$,
then we glue $Z$ and $G$ together over the boundaries, 
to obtain the rectangular lattice $W$.
Finally, we glue~$X$ and $W$ together over the boundaries, 
to obtain the rectangular lattice~$V$, 
which we call the \emph{triple gluing} of $G, Y,  Z, U$.

\begin{lemma}\label{L:gluecongfour}
A congruence $\bga$ of $V$ is uniquely associated with 
the four congruences $\bga_G$ of $G$,  $\bga_Y$ of $Y$, $\bga_Z$ of $Z$, 
and $\bga_U$ of $U$, 
satisfying the con\-dition that $\bga_G$ and $\bga_Y$ 
agree on the facing boundaries,
and the same for $\set{G,Z}$, $\set{Y,U}$, and $\set{Z,U}$.
\end{lemma}

This is a very easy lemma. For a formal proof of a more general statement, 
see the Appendix.

\section{Proof by Picture}\label{S:picture}

\begin{figure}[thb]
\centerline{\includegraphics{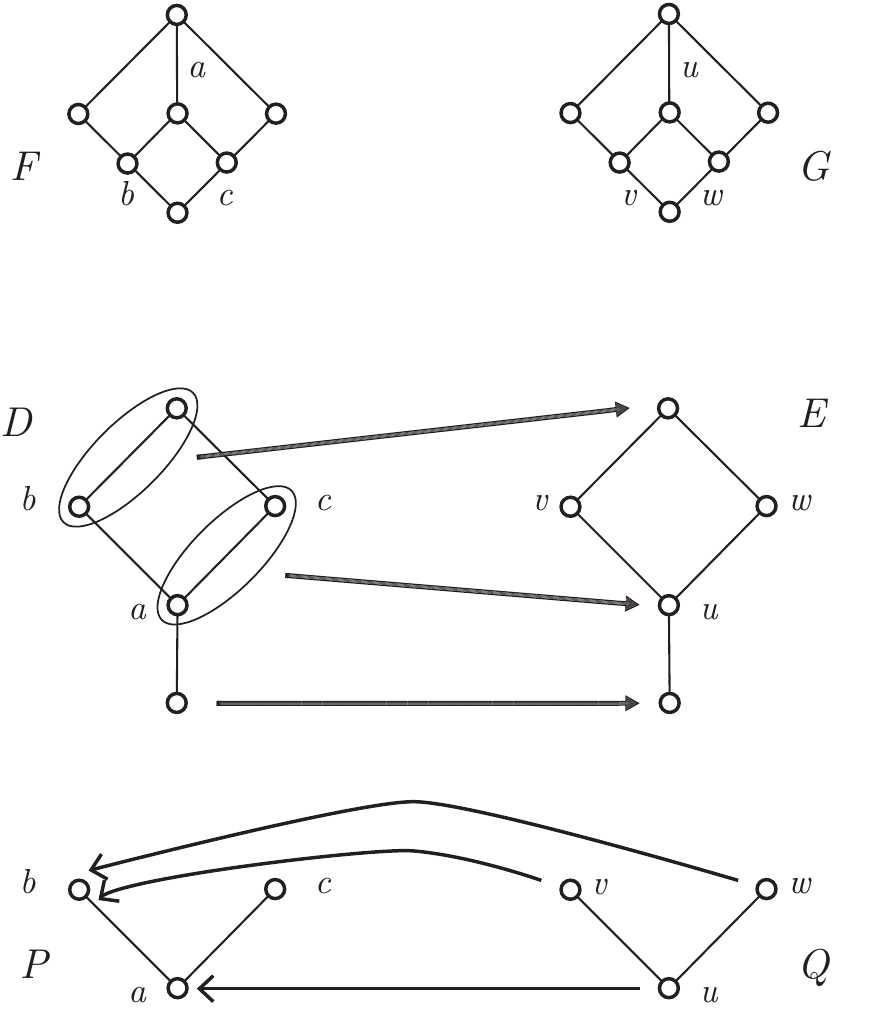}}
\caption{
(i) The rectangular lattices $F$ and $G$.\\
(ii) The distributive lattices $D = \Con F$, $E = \Con G$ 
and the bounded homomorphism $\gf \colon D \to E$.\\
(iii)  The ordered sets $P = \Ji D$, $Q = \Ji E$, 
and the isotone map $\Ji \gf \colon Q \to P$.
\label{F:setup}}
\end{figure}

\begin{figure}[p]
\centerline{\includegraphics[scale=.9]{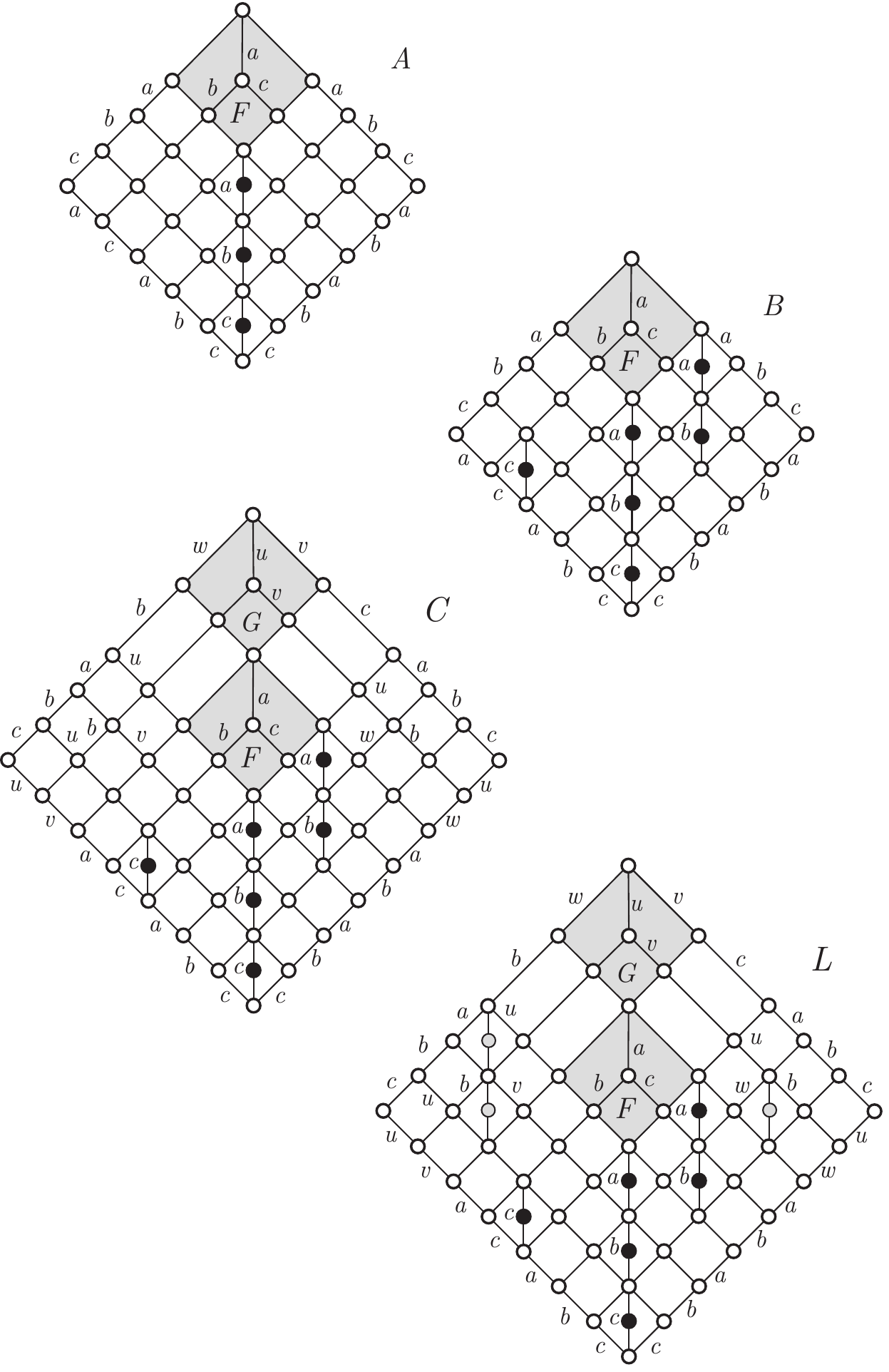}}
\caption{The rectangular lattices $A$, $B$,  $C$, and $L$.}
\label{F:all}
\end{figure}

For the ``Proof by Picture'', we choose 
\[
    F = G = \SfS 7.
\]
Let 
\begin{equation}
\begin{alignedat}{2}\label{F:DE}
    D &= \Con F, &E &= \Con G,\\ 
    P &= \Ji {D},\q  &Q &= \Ji{E},
\end{alignedat}
\end{equation}
and let $\gf$ the bounded homomorphism $D \to E$,
as in Figure~\ref{F:setup}(ii). 
Note that $\SfS 7$ is ``square'', 
as a result all the rectangular lattices look ``square'' in this 
and the subsequent diagrams of this chapter.

We construct the rectangular lattice $L$ of Theorem~\ref{T:Main} in a few steps,
as~illustrated in Figure~\ref{F:all}.

Step 1. We construct a congruence-preserving extension $B$ of $F$
in which all  join-irreducible congruences of $F$
appear on both upper boundaries of~$B$.

The distributive lattice $E = \Con G$ has $3$ join-irreducible congruences,
so we start with the glued sum of $3$ copies of $\three$ and $F$
and extend this to a rectangular lattice~$A$, see Figure~\ref{F:all}.

We add three more eyes to $A$ to form $\three$-s, 
to make sure that any two edges of the same color 
generate the same congruence; 
thus we obtain the rectangular lattice~$B$ of Figure~\ref{F:all}.

\vspace{-10pt}

In Figures~\ref{F:all}--\ref{F:second}, 
an $\three$, all whose edges are colored by $x$ are pictured as
 {\includegraphics[scale = .6]{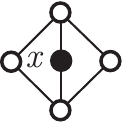}}.

Step 2. We form a rectangular lattice $C$ 
containing all the congruences of~$F$ and~$G$.
To do this, we form the glued sum of $B$ and $G$
and extend this to a rectangular lattice $C$ as in Figure~\ref{F:all}.

Step 3. We go back to Figure~\ref{F:setup}, 
to the element $v$ of $Q$.
The element~$b$ of $D$ is the smallest element mapped by $\gf$
to an element of $E$ with $\geq v$;
the element $b$ is join-irreducible, so it is in $P$.
We identify $b$ and $v$ 
(as per the discussion of the Birkhoff Representation Theorem),
by finding a cover preserving $\Csquare$ 
colored by $b$ and $v$ and adding an eye (colored gray in Figure~\ref{F:all}).
We proceed the same way with the element 
$w \in Q$, again finding $b \in E$, 
and identify $b$ and $w$. 
Finally, we do the same for the element $u  \in Q$ and
identify the congruences represented by $u$ and $a$.
Thereby, $L$ is a congruence-preserving extension of $F$, since each
join-irreducible congruence of $L$ is one of $a, b, c$.
This finishes the construction of the rectangular lattice~$L$,
as in Figure~\ref{F:all}.

Now we infer from BRT, 
that the map $\gf$ is represented by $\bga \restr F \to \bga \restr G$,
for $\bga \in \Con F$,
as required in Theorem~\ref{T:Main}.

\section{Proof of Theorem~\ref{T:Main}}\label{S:ProofMain}

We use the notation \eqref{F:DE}.

\subsection*{The first triple gluing}\label{S:first}

For the rectangular lattice $F$ of Theorem~\ref{T:Main}, 
let $\tup{bl}_F$, $\tup{br}_F$, $\tup{tl}_F$, and $\tup{tr}_F$ denote the number of elements 
of the bottom left and right, and top left and right boundaries, respectively,
and let $j$ denote the number of join-irreducible congruences.
If the lattice $F$ is understood, the subscripts may be omitted.

\begin{figure}[hbt]
\centerline{\includegraphics[scale=0.8]{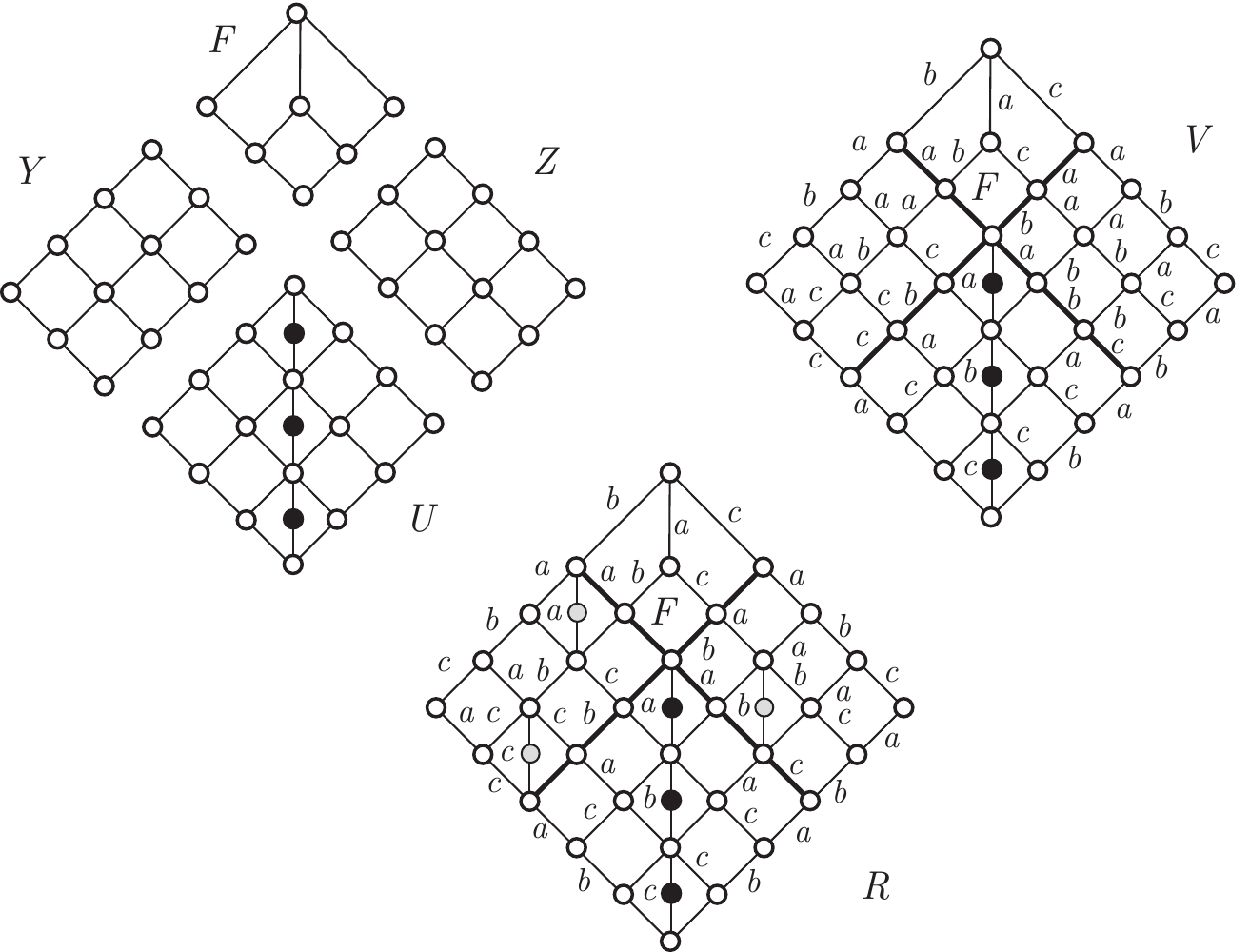}}
\caption{The first triple gluing}
\label{F:first}
\end{figure}

We prove the following statement.

\begin{lemma}\label{L:first}
The rectangular lattice $F$
has a congruence-preserving rectangular extension $R$
such that all join-irreducible congruences of $R$
appear on both upper boundaries of $R$.
\end{lemma}

\begin{proof}
Let us define the rectangular lattices, $Y$ and $Z$ as follows.
\[
   Y = \SC{\tup{bl}} \times \SC{j + 1},\q Z =  \SC{j + 1} \times \SC{\tup{br}}.
\]
We also need  the rectangular lattice $U$ that we obtain from  $\SC{j + 1}^2$
by adding eyes $j$-times 
to the covering $\SC 2^2$-s on the main diagonal.

We label~$F$; there are $j$ labels. 
We use these to label the $\three$-s on the main diagonal of $U$, bijectively.

Now we form the triple gluing of $F, Y, Z, U$ to form the rectangular lattice~$V$, 
as illustrated in Figure~\ref{F:first}. The labeling of $V$ is not a coloring, 
because two edges of the same label 
do not necessarily generate the same congruence.
Let $x$ be the label of a join-irreducible congruence of~$F$. 
There is an edge $A_x$ of label $x$ in the lower boundary of~$F$,
say, in the lower left boundary. 
Let $A_x'$ and  $B_x'$ be the edges of $Y$ 
that are identified with $A_x$ and~$B_x$, respectfully, in the triple gluing. 
Let $\bga_x$ and $\bgb_x$ denote the congruences of $Y$ 
generated by $A_x'$ and  $B_x'$, respectfully.

We claim that for the color $x$,
there is a covering square $S_x = \Csquare$ in $Y$,
that is colored by both $\bga_x$ and~$\bgb_x$. 
Indeed, take the trajectory~$\F r$ of~$Y$ containing~$A_x'$;
it~is a normal-up trajectory. 
Take the trajectory~$\F t$ of~$Y$ containing $B_x'$;
it~is a normal-down trajectory.
Therefore, $\F r$ and~$\F t$ intersect in a~covering square in $Y$
colored by both $\bga_x$ and~$\bgb_x$, as claimed.

We add an eye to $S_x$, as in Figure~\ref{F:first}. 
We do it for all colors~$x$ in~$F$, to~obtain the rectangular lattice $R$,
illustrated in Figure~\ref{F:first}.
Now the labelling is a coloring, the color $x$ of $F$ 
is the same as the color $x$ of~$U$.

Let $\bgg_l$ be the restriction of $\bga_x$ to $\SC{\tup{bl}_Y}$, 
the lower left boundary of $Y$. 
Similarly, associate with the congruence $x$ of $U$ the congruence of $Y$
generated by $C_x'$ that is identified with $C_x$, 
and restrict it to the lower right boundary of $Y$;
call it $\bgg_r$. 

Finally, $R$ is a congruence-preserving extension of $F$, as claimed in this lemma. 
To accomplish this, we~define the congruences $\bga_F$ of~$F$,  $\bga_Y$ of $Y$, 
$\bga_Z$ of $Z$, and $\bga_U$ of $U$, as follows: 

\begin{enumeratei}
\item The congruence $\bga_F$ of $F$ is $x$, equivalently, 
$\bga_F$ is generated in $F$ by the edge~$A_x$.
\item We define $\bga_Y = \bgg_l \times \bgg_r$.
\item The congruence $\bga_Z$ of $Z$ is defined symmetrically.
\item The congruence $\bga_U$ of $U$ is generated by the edge~$B_x \ci U$.
\end{enumeratei}
By construction, these congruences satisfy the conditions of Lemma~\ref{L:gluecongfour},
so there is a congruence $\bga_R$ extending all four, 
and this is the extension of the color $x$ to $R$.
This gives us that there is at least one such extension. 
The uniqueness follows from the fact the every edge of $R$ 
is perspective to an edge of $F$ or $U$.

We conclude by observing that the congruence extension property 
for join-irreducible congruences
is equivalent to the congruence extension property (for all congruences).
\end{proof}

\subsection*{The second triple gluing}\label{S:second}

We prove the following statement in this section.

\begin{lemma}\label{L:cep2}
The rectangular lattice $G$ has a congruence-preserving extension~$L$
such that $\id{0_G} = R$ \lp the rectangular lattice of Lemma~\ref{L:first}\,\rp.
\end{lemma}

\begin{figure}[htbp]
\centerline{\includegraphics[scale=.85]{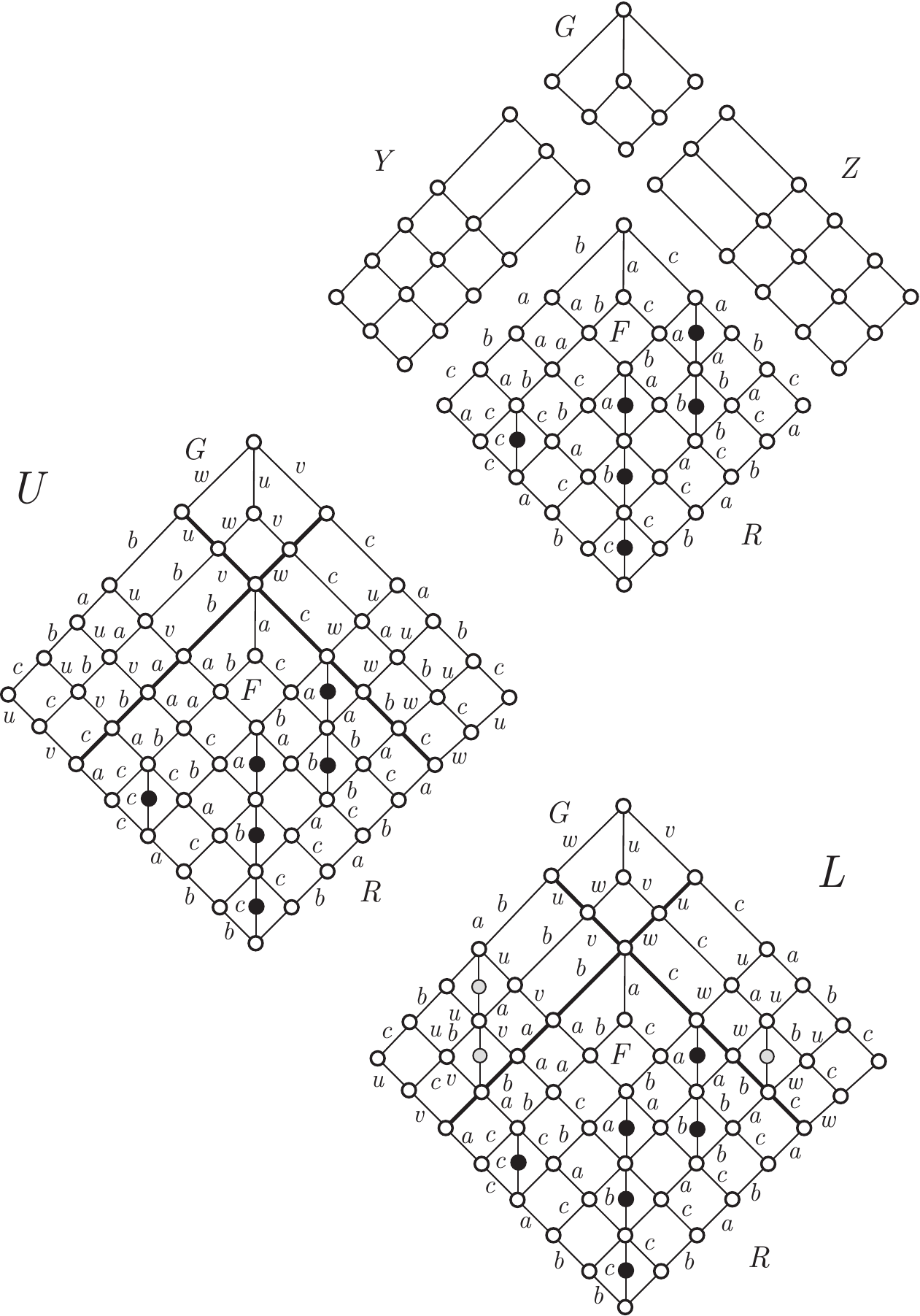}}
\caption{The second triple gluing.\\
(i) Top diagram: Setting up the second triple gluing.\\
(ii) Middle diagram: The result of the second triple gluing.\\
(iii) Bottom diagram: The rectangular lattice $L$.}
\label{F:second}
\end{figure}

\begin{proof}
The technical aspects of the proof are very similar to the proof of Lemma~\ref{L:first},
\emph{mutatis mutandis}.

We use three auxiliary rectangular lattices, defined as follows:
the rectangular lattice $R$ constructed in the previous section, and
\[
   Y = C_{\tup{bl}_G} \times C_{\tup{tl}_R},\q Z = C_{\tup{tr}_R} \times C_{\tup{br}_G}.
\]
Then we form the triple glued sum of $G, Y, Z, R$ 
to form the rectangular lattice~$U$, as illustrated in Figure~\ref{F:second}.(ii).

Recall that we use the notation \eqref{F:DE}.

Let $x$ be the color of a join-irreducible congruence of~$G$, that is, $x \in Q$. 
As~in Section~\ref{S:picture}, 
there is a smallest element $y \in D$ for which $\gf y \geq x$ holds,
namely, $\Ji \gf x \in P$.
There is an edge $A_x$ of color $x$ in the lower boundary of~$G$,
say, in the lower left boundary. 
There is also an edge $B_y$ of color $y$ in the upper left boundary of~$R$. 

Let $A_x'$ and  $B_y'$ be the edges of $Y$ 
that are identified with $A_x$ and~$B_y$, respectfully, in the second riple gluing. 
Let $\bga_x$ and $\bgb_y$ denote the congruences of $Y$ 
generated by $A_x'$ and  $B_y'$, respectfully.

 As the proof of Lemma~\ref{L:first},
we identify $x$ and $y$ 
by finding a cover preserving $\Csquare$  in~$Y$
colored by $x$ and $y$ and adding an eye (colored gray in Figure~\ref{F:all}). 
We do it for all colors~$x$ in $G$, to obtain the rectangular lattice~$L$,
illustrated in Figure~\ref{F:second}.

We verify the properties of $L$ 
as we verified the properties of $R$ in the previous section.
\end{proof}

\subsection*{Completing the proof of Theorem~\ref{T:Main}}

In the previous section, we constructed a rectangular lattice~$L$ 
containing the filter $G$ and the convex sublattice $F$ such that 
\ref{T:Main}(i) holds.
To complete the proof of Theorem~\ref{T:Main}, 
we have to verify that \ref{T:Main}(ii) also holds for $L$.
By BRT,
it is equivalent to verify that 
for $\bga \in Q$,  
the map $\Ji \gf$ is represented by the construction,
which is evident.

\newpage

\section{Ideal embeddabilty}\label{S:embeddabilty}
First some terminology. Let $G$ be a rectangular lattice.

We say $G$ is \emph{simple-embeddable}
if it is embeddable as an ideal in a simple rectangular lattice.

We say $G$ is \emph{abstractly ideal-representable} if, 
for any finite distributive lattice~$D$ 
and any bounded lattice homomorphism $\gf \colon D \to \Con G$, 
there is a rectangular lattice $L$ containing $G$ as an ideal, 
and an isomorphism
$\gh \colon D \to \Con L$ such that $\gf =\re \circ \gh$.

We say $G$ is \emph{concretely ideal-representable} if, for any rectangular lattice
$F$ and any bounded lattice homomorphism $\gf \colon \Con F \to \Con G$, there
is a rectagular congruence-preserving extension $L$ of $F$ that contains $G$ as
an idealand such that for any $\bga$ of $\Con L$,  
the map $\gf$ is $\bga\restr F \mapsto \bga\restr G$.

\begin{lemma}\label{L:meet}
Let $L$ be a lattice, let $I$ be an ideal in $L$, and let $\bga$ be a
meet-congruence on $I$. Extend $\bga$ to an equivalence relation $\bgb$
on $L$ by setting the  equivalence class of any $x \in L \setminus I$ to be the
singleton $\set{x}$. Then $\bgb$ is a meet-congruence on $L$
\end{lemma}
\begin{proof}
We need only show that for $x$, $y$, $z \in L$, if $x \neq y$ and
$\cng x = y (\bgb)$, then $\cng x \mm z = y \mm z (\bgb)$. But
then $x$, $y \in I$, and $\cng x = y (\bga)$. Now,
\[
   z' = (x \jj y) \mm z \in I
\]
and
\[
   x \mm z = x \mm z', \quad y \mm z = y \mm z'.
\]
Thus $\cng x \mm z = y \mm z (\bga)$, that is, $\cng x \mm z = y \mm z (\bgb)$
\end{proof}

\begin{lemma}\label{L:ideal}
Let $L$ be a rectangular lattice, and let the ideal $I$ of $L$ also be a
rectangular lattice. Then $\lc(I)$ is
in the left lower chain of $L$ and $\rc(I)$ is in the
right lower chain of $L$.
\end{lemma}
\begin{proof}
The corners of $I$, the elements $\lc(I)$ and $\rc(I)$, are join-irreducible in the ideal $I$, and so
are join-ireducible elements of $L$. But any join-irreducible element of $L$ is either an
eye or an element of one of the two lower chains of $L$. Being incomparable, they
cannot both lie in the same chain.
 
On the other hand, we claim that neither can be an eye in $L$. An eye $e$ in $L$
has a unique upper cover in $L$ which has two additional lower covers, one to the left of
$e$ and one to the right of $e$. But each of $\lc(I)$ and $\rc(I)$ has an upper
cover in one of the upper chains of the ideal $I$, and so, by the uniqueness
of the upper cover, if $\lc(I)$ were an eye, we would get the contradiction
that $I$ contains an element to the left of $\lc(I)$, and, similarly for $\rc(I)$,
an element to the right of $\rc(I)$.

Thus $\lc(I)$ is in the left lower chain of $L$ and $\rc(I)$ 
is in the right lower chain of~$L$.
\end{proof}

We note the following.

\begin{lemma}\label{L:join}
Let $L$ be a rectangular lattice and let $x \in L$ not be an eye. Then
$x = (x \mm \lc(L)) \jj (x \mm \rc(L)$.
\end{lemma} 

\begin{lemma}\label{L:diff}
Let $L$ be a rectangular lattice, and let the ideal $I$ of $L$ also be rectangular.
Then for each $x \in L \setminus I$, either $x > \lc(I)$ or $x > \rc(I)$.
\end{lemma}
\begin{proof}
We first consider the case where $x$ is not an eye. Then, by Lemma~\ref{L:join},
$x = (x \mm \lc(L)) \jj (x \mm \rc(L)$. Now, $\Cll(L)$  is a chain containing both
$x \mm \lc(L)$ and, by Lemma~\ref{L:ideal}, $\lc(I)$. Thus
\[
   \text{either } x \mm \lc(L)  \leq  \lc(I)\quad\text{ or }\quad x \mm \lc{L}  >  \lc(I).
\]
Similarly,
\[
   \text{either } x \mm \rc(L)  \leq  \rc(I)\quad\text{ or }\quad x \mm \rc{L}  >  \rc(I).
\]
But, if $x \mm \lc(L)  \leq  \lc(I)$ and $x \mm \rc(L)  \leq  \rc(I)$, we get the contradiction
$x \leq \lc(I) \jj \rc(I) \in I$. Thus, either $x \mm \lc{L}  >  \lc(I)$, and so $x > \lc(I)$,
or $x \mm \rc{L}  >  \rc(I)$, and so $x > \rc(I)$, establishing our claim if $x$ is not
an eye.

On the other hand, if $x$ is an eye, then $x$ has a unique lower cover $x'$ in $L$,
and $x'$ is not an eye. If $x'  \in \Cul{I} \uu \Cur{I}$, then
\[
   \text{either } x > x' \geq \lc(I) \quad\text{ or }\quad x > x' \geq \rc(I).
\]
Otherwise, $x' \in L \setminus I$, and so, by the first case of the proof,
either $x ' > \lc(I)$ or $x' >  \rc(I)$, and thus either $x > \lc(I)$ or
$x > \rc(I)$, establishing our claim also when $x$ is an eye.
\end{proof}

\begin{lemma}\label{L:nosimple}
Let $L$ be a rectangular lattice, let $I$ be an ideal in $L$ which is also rectangular
 and let $\bga$ be a congruence on $I$ that collapses no interval in either upper chain 
of  $I$. Extend $\bga$ to an equivalence relation $\bgb$
on $L$ by setting the  equivalence class of any $x \in L \setminus I$ to be the
singleton $\set{x}$. Then $\bgb$ is a congruence on $L$
\end{lemma}

\begin{proof}
By Lemma~\ref{L:meet}, the relation $\bgb$ is preserved by meet. So, we need only
show that it is preseved by join. We need only consider elements $x$, $y$, $z \in L$
with $x < y$ and $\cng x = y (\bgb)$. Being distinct, $x$, $y \in I$, by definition of
$\bgb$, and so $\cng x = y (\bga)$. If $z \in I$ then $\cng x \jj z = y \jj z (\bga)$
and we are done.

Otherwise, by Lemma~\ref{L:diff}, we may assume that $z > \lc(I)$. But then
\[
x \jj z = x \jj \lc(I) \jj z \quad\text{and}\quad  y \jj z = y \jj \lc(I) \jj z.
\]
Now $ x \jj \lc(I)$ and $y \jj \lc(I)$ are in the left upper chain of $I$ and 
$\cng x \jj \lc(I) = y \jj \lc(I) (\bga)$. By our  condition on $\bga$ we get
$x \jj \lc(I) = y \jj \lc(I)$, and so $x \jj z = y \jj z$. Thus $\bgb$ is also preserved
by join.
\end{proof}

We have an immediate corollary:
\begin{lemma}\label{L:noideal}
Let the rectangular lattice $G$ have a congruence that collapses no interval
in either upper chain. Then $G$ cannot be embedded as a ideal in any simple
rectangular lattice.
\end{lemma}

\begin{lemma}\label{L:concrete}
Let $G$ be a rectangular lattice such that each non-trivial congruence collapses some
edge in one of its upper chains.

Let $F$ also be a rectangular lattice
and let $\gf \colon \Con{F} \to \Con{G}$
be a~bounded lattice homomorphism. 
Then there is a rectangular lattice~$L$ with $G$ as an ideal
and $F$ as a filter satisfying the following two conditions\co
\begin{enumeratei}
\item the lattice $L$ is a congruence-preserving extension of $F$\tup{;}
\item for $\bga$ of $\Con L$,  
the map $\gf$ is $\bga\restr F \mapsto \bga\restr G$.
\end{enumeratei}
\end{lemma}
\begin{proof}
(Outline) We proceed as in the proof of Theorem~\ref{T:Main}, 
\emph{mutatis mutandis}
We~have $\Ji{\gf} \colon \Ji{G} \to \Ji{F}$. Embed $F$ as a filter in a
preserving rectangullar extension~ $F'$ where each non-trivial congruence
appears as a color on each lower chain. We then glue $F'$ to the top of
$G$ and add the two "flaps". Finally, insert eyes in the flaps so that, for
each join-irreducible congruence $\bga$ of~$G$, each edge on either upper
chain of~$G$ that is colored by $\bga$ is connected to an edge on the
corresponding lower chain of $F'$ that is colored by $\Ji{\gf}(\bga)$. The
lattice $L$ is the required rectangular lattice.
\end{proof}
 By Lemmas~\ref{L:noideal} and \ref{L:concrete}, we have our main theorem.
\begin{theorem}
Let $G$ be a rectangular lattice. The following conditions are equivalent\co
\begin{enumeratei}
\item each non-trivial congruence of $G$ collapses some edge in one of
its upper chains\tup{;}
\item $G$ is concretely ideal-representable\tup{;}
\item $G$ is abstractly ideal-representable\tup{;}
\item $G$ is simple-embeddable.
\end{enumeratei}
\end{theorem}

\section{Discussion}\label{S:Discussion}

\subsection{Rectangular extensions}\label{S:rectangular}

In Section~\ref{S:picture}, Proof by Picture, 
we write: ``extend this to a rectangular lattice'' in Step 1.
The is unambiguous but not very precise. 
The formal way is to use triple gluing.

G. Gr\"atzer and E.~Knapp~\cite{GKn09} 
prove the existence of a rectangular congruence-preserving extension
of an SPS lattice. 
The uniqueness of a rectangular extension 
is discussed in G.~Cz\'edl~\cite{gC17}.

\subsection{$C_1$-diagrams}\label{S:diagrams}

G. Cz\'edli \cite{gC17} introduced $C_1$-diagrams for slim rectangular lattices.
Informally, a diagram is a  $C_1$-diagram, 
if all edges are drawn at directions of $(1, 1)$ and $(1, -1)$,
\emph{normal edges}, 
except the middle edges of $\seven$-s are drawn at directions 
$(\cos \alpha, \sin \alpha)$ with $\pi/2 < \alpha <3\pi/2$,  \emph{steep edges}.
G.~Cz\'edli proves that all slim rectangular lattices have  $C_1$-diagrams,
see  \cite{gC17}.
See also G.~Gr\"atzer~\cite{gG20}.

All the diagrams of rectangular lattices of this paper are $C_1$-diagrams.

\subsection{Automorphisms}\label{S:Automorphisms}

It is easy to see that the lattice $L$ of Theorem~\ref{T:Main}
can be constructed with a given automorphism group.

\section*{Appendix}

Any gluing can be described as a lattice $L$ with an ideal $A$ and
a filter $B$ such that $A \ii B \neq \es$; the lattice is the gluing of any lattice $A'$ isomorphic to $A$
and any lattice $B'$ isomorphic to $B$ over the filter of $A'$ corresponding to $A \ii B$ and
ideal of $B'$ isomorphic to $A \ii B$. All that is required is that $A$ and $B$ have a nonempty intersection.

We recall the following characterization of how congruences on $A$, $B$, and $L$ relate.

\begin{lemma}\label{L:gluing}
Let the lattice $L$ be a gluing of an ideal $A$ and a filter $B$. Let $\bga_A$ be a congruence
of $A$ and $\bga_B$ a congruence of $B$. Then there is a congruence $\bga$ of $L$ with
$\bga_A = \bga \restr A$ and $\bga_B = \bga \restr B$ if{f}
\[
\bga \restr (A \ii B) = \bga_B \restr (A \ii B),
\]
and, in that event, 
\[
\bga = \bga_A \uu \bga_B \uu (\bga_A \circ \bga_B) \uu (\bga_B \circ \bga_A),
\]
uniquely determined by $\bga_A$ and $\bga_B$.
\end{lemma}

We now introduce triple gluing.

Let $L$ be a lattice and let $c \in L$. We assume that $L$ is the union of four convex sublattices,
$U$, $V$, the ideal $X = \Dg{c}$, and the filter $Y = \Ug{c}$.

We assume that the following properties hold:
\begin{enumeratei}
\item[(Pi)] $A = U \uu X$ is an ideal of $L$.
\item[(Pii)] $B = V \uu Y$ is a filter of $L$.
\item[(Piii)] $U$ is a filter of $A$ and $X$ is an ideal of $A$.
\item[(Piv)] $V$ is an ideal of $B$ and $Y$ is a filter of  $B$.
\item[(Pv)] $U \ii V = \set{c}$.
\end{enumeratei}

\begin{figure}[htb]
\centerline{\includegraphics{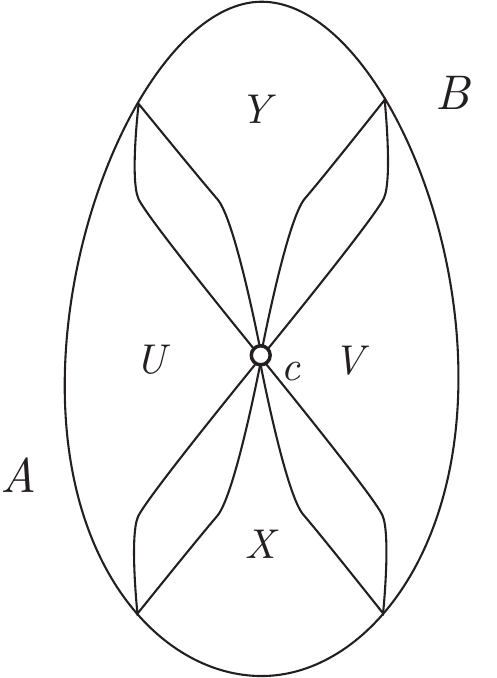}}
\caption{Triple gluing}\label{F:Htriple}
\end{figure}

By (Pv), $c$ is an element of each of $X$, $Y$, $U$, and $V$; thus no pairwise intersection
is empty. Then
(Pi) and (Pii) state that $L$ is a gluing of $A$ and $B$ over $A \ii B$. 
(Piii) states that~$A$ is a gluing of $X$ and $U$ over $X \ii U$, 
and (Piv) states that $B$ is a gluing of $V$ and $Y$ over $V \ii Y$. 
We call such a configuration a \emph{triple gluing}; see Figure~\ref{F:Htriple}.

Property (Pv) will be utilized in the following results.

\begin{lemma}\label{L:UY}
$U \uu Y$ is a sublattice of $L$, with $U$ an ideal and $Y$ a filter, so
$U \uu Y$ is a gluing of $U$ and $Y$ over $U \ii Y$.

Dually, $X \uu V$ is a sublattice of $L$ with $V$ a filter and $X$ an ideal, 
so $X \uu V$ is a gluing of $V$ and $X$ over $V \ii X$.
\end{lemma}

\begin{proof}
By duality, we need only prove the first statement.

We first show that $U \uu Y$ is a sublattice of $L$. Let $x, y \in U \uu Y$.
If $x, y \in U$ (respectively, $x, y \in Y$), 
then $x \mm y, x\jj y \in U$ (respectively $\in Y$),
since $U$ and $V$ are sublattices of $L$.

Otherwise, we may assume that $x \in U$ and $y \in Y$. Then $x \jj y \in Y$, since $Y$ is a
filter of $L$. So, we need only show that $x \mm y \in U$. We have $y \geq c$, by definition of~$Y$.
Thus,
\[
   x \mm c \leq x \mm y \leq x.
\]
Then $c \in U$ follows from (Pv). Thus $x \mm c, x \in U$, and, since $U$ is a
convex sublattice of $L$, we conclude that $x \mm y \in U$. Thus $U \uu Y$ is a sublattice
of $L$.

Now $Y$ was defined as a filter of $L$, and so is clearly a filter of $U \uu Y$.

We then need only show that $U$ is an ideal of $U \uu Y$. By (Pi),
$(U \uu Y) \ii A$ is an ideal of $U \uu Y$. But
\[
   (U \uu Y) \ii A = U \uu (Y \ii X) = U \uu \set{c} =U,
\]
since $Y \ii X = \set{c}$ by definition of $X$ and $Y$, and since $c \in U$. Thus,
$U$ is indeed an ideal of the sublattice $U \uu Y$.
\end{proof}

\begin{lemma}\label{L:main}
Let the lattice $L$ be the triple gluing of $U$, $V$, $X$, $Y$. 
Let $\bga_X$ be a congruence of $X$, 
$\bga_Y$ a congruence of $Y$, $\bga_U$ a congruence of $U$, 
and $\bga_V$ a congruence of $V$. Furthermore, let
\begin{align}
   \bga_X \restr (U \ii X) &= \bga_U \restr (U \ii X),\label{E:UX} \\
   \bga_X \restr (V \ii X) &= \bga_V \restr (V \ii X),\label{E:VX}\\
   \bga_Y \restr (V \ii Y) &= \bga_V \restr (V \ii Y),\label{E:VY}\\
   \bga_Y \restr (U \ii Y) &= \bga_U \restr (U \ii Y).\label{E:UY}
\end{align}
Then there exists a unique congruence $\bga$ on $L$ such that $\bga \restr X = \bga_X$,
$\bga \restr Y = \bga_Y$, $\bga \restr U = \bga_U$, and $\bga \restr V = \bga_V$.
\end{lemma}
\begin{proof}
We apply Lemma~\ref{L:gluing} successively  to the lattice $A = U \uu X$,
the lattice $B = V \uu Y$, and the lattice $L = A \uu B$.

By \eqref{E:UX}, $\bga_X$, and $\bga_U$ extend uniquely to the congruence
\begin{equation}\label{E:beta}
\bgb = \bga_X \uu \bga_U \uu(\bga_X \circ \bga_U) \uu (\bga_U \circ \bga_X)
\end{equation}
on
$A $. Similarly, by \eqref{E:VY}, 
the congruences  $\bga_Y$ and $\bga_V$ extend
uniquely to the congruence
\begin{equation*}
   \bgg = \bga_Y \uu \bga_V \uu (\bga_Y \circ \bga_V) \uu (\bga_V \circ \bga_Y)
\end{equation*}
on $B$.

To conclude the proof, we have to show that 
$\bgb$ and $\bgg$ extend uniquely to a congruence $\bga$ of $L$, 
that is, that
\begin{equation*}
   \bgb \restr (A \ii B) = \bgg \restr (A \ii B).
\end{equation*}
By duality, it suffices to show that
\begin{equation}\label{E:AB}
\bgb \restr (A \ii B) \leq \bgg \restr (A \ii B),
\end{equation}
that is, that for any pair $x < y$ of elements of $A \ii B$ where
$\cng x = y (\bgb)$,
it follows that $\cng x = y (\bgg)$. We first determine $A \ii B$. By definition of
$X$ and $Y$, $X \ii Y = \set{c}$. Then, by (Pv), we get
\[
   A \ii B = (V \ii X) \uu (U \ii Y).
\]

If $x, y \in V \ii X$, then $\cng x = y (\bga_X)$, and by \eqref{E:VX},
$\cng x = y (\bga_V)$ and so $\cng x = y (\bgg)$.

Similarly, if  $x, y \in U \ii Y$,  we conclude, by \eqref{E:UY}, that  $\cng x = y (\bgg)$.

Otherwise, since $X = \Dg{c}$  and $Y=\Ug{c}$, we conclude that
$x \in V \ii X$ and $y \in U \ii Y$, and so that $x < c <y$. Then,
$\cng c = y (\bgb)$. But, $c, y \in U \ii Y$ and
\[
   \bgb \restr (U \ii Y) = \bga_U \restr (U \ii Y) = \bga_Y \restr (U \ii Y) \ci \bgg.
\]
Thus
$
   \cng c = y (\bgg).
$

Furthermore, by \eqref{E:beta}, $\cng x = y (\bga_X \circ \bga_U)$, that is,
there is $z \in X \ii U$ with
$\cng x = z (\bga_X)$  and $\cng z = y(\bga_U)$.

Now $x \in V \ii X$, a filter in $X$, and $z \in U \ii X$, another filter in $X$, Thus
$x \jj z \in (U \ii X) \ii (V \ii X) = \set{c}$, that is $x \jj z = c$.  Then
$
   \cng x = c (\bga_X).
$
Since $x, c \in V \ii X$, $\cng x = c (\bga_V)$ by \eqref{E:VX}. Thus
$
   \cng x = c (\bgg).
$
Therefore
$
   \cng x = y (\bgg),
$
thereby concluding the proof of \eqref{E:AB}, and so concluding the proof of the lemma.
\end{proof}

\end{document}